\documentclass{article}
\usepackage{cite}
\usepackage[pdftex]{graphicx}
\usepackage{amsmath}	
\usepackage{array}

\usepackage{stfloats}
\usepackage{url}
\usepackage{color}
\usepackage[ruled,linesnumbered]{algorithm2e}
\usepackage{multirow}
\usepackage{amsthm}
\usepackage{amssymb}
\usepackage{amsmath}
\usepackage{mathrsfs}

\newcommand{\isr}{\mathtt{ISR}}

\newcommand{\crb}{\mathtt{CRLB}}

\newcommand{\cov}{\mathtt{cov}}

\newcommand{\A}{{\bf A}}
\newcommand{\W}{{\bf W}}

\newcommand{\y}{{\bf y}}
\newcommand{\tr}{\mathtt{tr}}
\renewcommand{\a}{{\bf a}}
\newcommand{\z}{{\bf z}}
\newcommand{\x}{{\bf x}}

\newcommand{\g}{{\bf g}}
\newcommand{\h}{{\bf h}}
\newcommand{\I}{{\bf I}}

\newcommand{\C}{{\bf C}}

\newcommand{\bkappa}{\boldsymbol{\kappa}}
\newcommand{\okappa}{\overline{\kappa}}
\newcommand{\bnabla}{\boldsymbol{\nabla}}
\newcommand\bigzero{\makebox(0,0){\text{\huge0}}}

\newcommand{\p}{\partial}

\title{Cram\'er-Rao Bounds for Complex-Valued Independent Component Extraction: Determined and Piecewise Determined Mixing Models}
\author{
V\'aclav Kautsk\'y$^{1,2}$, Zbyn\v{e}k Koldovsk\'{y}$^{2}$, Petr 
Tichavsk\'{y}$^3$, and Vicente Zarzoso$^4$}

\begin{document}

\maketitle

\footnotetext{This work was supported by The Czech Science Foundation through 
Project No.~17-00902S and by 
the United States Department of the Navy, Office of Naval Research Global, 
through Project No.~N62909-18-1-2040.
	\\$^1$Faculty of Nuclear Sciences and Physical Engineering,
	Czech Technical University in Prague, e-mail: {\tt kautsvac@fjfi.cvut.cz}
	\\$^2$Faculty of Mechatronics, Informatics, and Interdisciplinary Studies, 
	Technical University of Liberec,
	Studentsk\'a 2, 461 17 Liberec, Czech Republic
	\\$^3$Institute of Information Theory and Automation,
	P.O.Box 18, 182 08 Prague 8, Czech Republic
	\\$^4$GEII Department and the I3S Laboratory, University of Nice Sophia Antipolis, France.}

\begin{abstract}
This paper presents Cram\' er-Rao Lower Bound (CRLB) for the complex-valued Blind Source 
Extraction (BSE) problem based on the assumption that the target signal is independent of  
the other signals. Two instantaneous mixing models are considered. First, we consider the standard determined mixing model used in Independent Component Analysis (ICA) where the mixing matrix is square and non-singular and the number of the latent sources is the same as that of the observed signals. The CRLB for Independent Component Extraction (ICE) where the mixing matrix is re-parameterized in order to extract only one independent target source is computed. The target source is assumed to be non-Gaussian or non-circular Gaussian while the other signals (background) are circular Gaussian or non-Gaussian. The results confirm some previous observations known for the real domain and bring new results for the complex domain. Also, the CRLB for ICE is shown to coincide with that for ICA when the non-Gaussianity of background is taken into account. %unless the assumed sources' distributions are misspecified. 
Second, we extend the CRLB analysis to piecewise determined mixing models. Here, the observed signals are assumed to obey the determined mixing model within short blocks where the mixing matrices can be varying from block to block. However, either the mixing vector or the separating vector corresponding to the target source is assumed to be constant across the blocks. The CRLBs for the parameters of these models bring new performance bounds for the BSE problem.

%A recently considered re-parametrization of the mixing model in Independent Component Extraction (ICE) is applied used for the standard determined mixing model  where, compared to Independent Component Analysis (ICA), the mixing model contains minimum number of parameters needed for the extraction problem.  while the other signals, which are not separated from each other, are modeled as  circular Gaussian or non-Gaussian. A CRLB-induced Bound (CRIB) for Interference-to-Signal Ratio (ISR) is derived and compared with similar bound for ICA. Numerical simulations with selected BSE algorithms show the correspondence between empirical results and the theory.
\end{abstract}

\section{Introduction}
\subsection{Problem Statement}
In Blind Source Separation (BSS), the instantaneous linear mixing model 
\begin{equation}\label{ICA}
\x = \A {\bf u}
\end{equation} 
is studied, where $\x$ is a $d\times 1$ vector representing $d$ observed signals, ${\bf u}$ is a $n\times 1$ vector of source signals, and $\A$ is a $d\times n$ mixing matrix. The goal of BSS is to  separate ${\bf u}$ from $\x$ using only information provided by the observed samples \cite{comon2010handbook}. Blind Source Extraction (BSE) aims at extracting only one source referred to as source of interest (SOI), while the other signals in $\x$ are called background. 
In this paper, complex-valued signals and parameters will be considered.

Independent Component Analysis (ICA) is a popular BSS method based on the assumption that the source signals are mutually independent. The $j$th source signal, $j=1,\dots,n$, $u_j$ (the $j$th element of ${\bf u}$) is modeled as a random variable with the probability density function (pdf) $p_j(\cdot)$, and the observed samples of $\x$ are assumed to be identically and independently distributed.
In the standard model, the determined case is considered where the number of sources is the same as that of the observed signals, $n=d$, and $\A$ is square $d\times d$ non-singular matrix. Here, the estimation of $\A$ and of $\A^{-1}$ is equivalent with the separation of ${\bf u}$, which is done through finding a square de-mixing matrix ${\bf W}$ such that ${\bf y}={\bf Wx}$ are as independent as possible. The identifiability and separability conditions were analyzed in \cite{complexErikson}. 

In this paper, we focus on the BSE problem where the SOI should be extracted based on the assumption of its independence from the background, a problem closely related to ICA. We compute Cram\'er-Rao Lower Bounds (CRLB) in order to analyze performance limitations of three mixing models, two of which were only recently considered in the literature \cite{koldovsky2019icassp}. The resulting bounds are compared between each other and also with the similar bound for the standard ICA. 

The paper has two parts. In the first part, the standard determined mixing scenario is considered, where the BSE problem is formulated through the recently proposed approach called Independent Component Extraction (ICE) \cite{koldovsky2018TSP}. Here, a particular parametrization of the mixing system is considered, which is designed for extracting only the first source $u_1$ from \eqref{ICA} playing the role of the SOI (without any loss on generality). 
Specifically, the mixing matrix and its inverse (de-mixing) matrix are parameterized, respectively, as
\begin{equation}\label{mixing}
\A_{\rm ICE}=\begin{pmatrix}{\bf a} & {\bf Q}\end{pmatrix}=
\begin{pmatrix}
\gamma & {\bf h}^H\\
\g & \frac{1}{\gamma}\left({\bf g}{\bf h}^H-{\bf I}_{d-1}\right)
\end{pmatrix}
\end{equation}
\begin{equation}\label{demixing}
\W_{\rm ICE}=\A_{\rm ICE}^{-1}=
\begin{pmatrix}
{\bf w}^H\\
{\bf B}
\end{pmatrix}=
\begin{pmatrix}
{\beta}^* & {\bf h}^H\\
\g & -\gamma{\bf I}_{d-1}
\end{pmatrix},
\end{equation}
where ${\bf a}$ denotes the first column of $\A$, which is the mixing vector related to $u_1$ partitioned as  ${\bf a}=[\gamma;{\bf g}]$, and ${\bf w}$ is the separating vector such that ${\bf w}^H\x=u_1$, partitioned as ${\bf w}=[\beta;{\bf h}]$. ${\bf I}_d$ denotes the $d\times d$ identity matrix, and $\beta$ and $\gamma$ are linked through 
\begin{equation}\label{betagamma}
{\beta}^*\gamma=1-{\bf h}^H{\bf g}.
\end{equation} 
This parametrization does not mean any restriction in the sense that $\A$ from \eqref{ICA} must obey the structure given by \eqref{mixing} in order to extract $u_1$. In fact, the extraction of the background subspace is ambiguous (any transformation of that subspace does not influence the independence of the background from the SOI), so \eqref{mixing} resp. \eqref{demixing} is just a particular choice that guarantees that ${\bf Ba}={\bf 0}$. % Algorithms derived in \cite{koldovsky2018TSP} through ICE (assuming that the background is circular Gaussian) were shown to be closely related to, for example, well-known One-Unit FastICA \cite{hyvarinen1999}, which does not assume any particular structure of $\A$. 
The ICE formulation enables us to compute the CRLB as we did in \cite{kautsky2017} for the real-valued case and Gaussian background. The contribution here compared to \cite{kautsky2017} is that the bound is derived for the complex-valued case and it involves also the non-Gaussian background.

In the second part of this paper, we compute the CRLBs for two {\em piecewise determined mixing models} that are designed for dynamic mixtures. Here, it is assumed that the observed samples of mixed signals can be partitioned into $M$ blocks where the samples in each block obey the standard determined model \eqref{ICA}. The $m$th block, $m=1,\dots,M$, is thus described by
\begin{equation}
\x^{m} = \A^{m}  {\bf u}^{m},
\label{newmodel}
\end{equation}
where the source signals ${\bf u}^m=[u_1^m,\dots,u_d^m]^T$ are independent. The mixing matrices $\A^{1},\dots,\A^{M}$ as well as the source signals (their distributions) can be varying from block to block\footnote{The equation \eqref{newmodel} is formally identical with the mixing model studied in Frequency-domain ICA \cite{smaragdis1998}, Independent Vector Analysis \cite{kim2007} or in joint BSS. There, the problem of joint blind separation of a set of instantaneous mixtures is considered, and $m$ plays the role of the mixture (dataset) index (e.g. the frequency bin index).}. The model thus involves dynamic mixing as well as a special underdetermined case (more sources than sensors) since there can be up to $M\times d$ sources. The fact that the mixtures are determined within the blocks brings the advantage of tractability of the analytic computation of the CRLB. % 

%The motivation for this model resides in the fact that the entire mixture (all blocks) can involve up to $M\times d$ sources. The model thus constitutes a special underdetermined\footnote{Difficult real-world signal mixtures are underdetermined, that is, they involve more sources than sensors ($n>d$).} case provided that $M>1$ (and coincides with \eqref{ICA} when $M=1$). The sub-models of blocks are invertible, which brings the advantage of tractability of the analytic computation of the CRLB. % and of developing efficient algorithms; see \cite{koldovsky2019icassp}.

However, without any further assumption, \eqref{newmodel} corresponds to a sequential application of the standard mixing model, which is straightforward for on-line signal processing but does not bring any advantage. Therefore, we propose special parametrizations useful for the BSE problem assuming that the SOI is active in all blocks and some mixing parameters related to the SOI are joint to all the blocks. Specifically, we parametrize $\A^1,\dots,\A^M$ similarly to \eqref{mixing} and consider two special variants: 
\begin{align}
\A^m_{\rm CMV}&=
\begin{pmatrix}
\gamma & ({\bf h}^m)^H\\
\g & \frac{1}{\gamma}\left({\bf g}({\bf h}^m)^H-{\bf I}_{d-1}\right)
\end{pmatrix}, \label{CMV} \\
\A^m_{\rm CSV}&=
\begin{pmatrix}
\gamma^m & {\bf h}^H\\
\g^m & \frac{1}{\gamma}\left({\bf g}^m{\bf h}^H-{\bf I}_{d-1}\right).
\end{pmatrix}. \label{CSV}
\end{align}
The models will be referred to as Constant Mixing Vector (CMV) and Constant Separating Vector (CSV), respectively, because, in CMV, the mixing vectors ${\bf a}^1,\dots,{\bf a}^M$ are constant over blocks and are equal to ${\bf a}$, and, in CSV, the separating vectors ${\bf w}^1,\dots,{\bf w}^M$ are all equal to ${\bf w}$. CMV is useful for situations where the SOI is a static source while the background is varying. CSV involves a moving SOI (varying mixing vector) under the assumption that 
a constant separating vector such that extracts the signal from all blocks exists. These models have been considered for the first time in \cite{koldovsky2019icassp}, where they were applied to blind audio source extraction. This paper provides their theoretical analysis through the CRLB theory.

%should provide a spatial filter that attenuates all directions where the source did not appear\footnote{Such filter could be seen as a blind estimator of the Linearly Constrained Minimum Variance (LCMV) beamformer \cite{vantrees2002}.}; see also \cite{koldovsky2019icassp}.

\subsection{State-of-the-Art}
\subsubsection{Independence-based BSS/BSE methods}
BSE methods based on signals' non-Gaussianity had been studied even before ICA was formulated \cite{huber1985,herault1987} in the Comon's pioneering paper \cite{comon1994}. Then, the theory of ICA has been established since 90s; see, e.g., \cite{cardoso1998,hyvarinen2001,cichocki2002,comon2010handbook}. The relation of the non-Gaussianity based BSE methods has been described through information theory and the properties of the Kullback-Leibler divergence (mutual information) and entropy \cite{cover}. ICE is a recent revision of this relation based on the algebraic mixing model \eqref{mixing} and maximum likelihood estimation \cite{koldovsky2018TSP}.

ICA has been used for blind separation of convolutive mixtures in the frequency domain \cite{smaragdis1998}, where the mixture is transformed into a set of complex-valued instantaneous mixtures (one mixture per frequency). The problem, called Frequency-Domain ICA (FDICA), is formally described by \eqref{newmodel}, however, $m$ plays the role of the frequency bin index. When ICA is applied separately to each mixture, the indeterminacy of the order of separated component gives rise to the permutation problem \cite{sawada2004sap} (the separated frequency components must be reordered in order to separate the signals in the frequency domain). 

To avoid the permutation problem, Independent Vector Analysis (IVA) has been proposed \cite{kim2007}. Here, the algebraic model remains the same as in FDICA while the statistical model involves the assumption that independent components belonging to the same source are  mutually dependent and form so-called vector components. The idea of IVA have become very popular due to its wide applicability far beyond audio source separation \cite{adali2015,chen2016}. Its variant for BSE (Independent Vector Extraction - IVE) appeared, e.g., in \cite{lee2007fast} and has been recently formulated in \cite{koldovsky2018TSP}. 

Another recent advancement in this line represents Independent Low Rank Matrix Analysis (ILRMA) where the statistical model of a vector component (representing one source) assumes that its spectrogram has a low-rank structure. For example, ILRMA combines IVA and Nonnegative Matrix Factorization (NMF) in \cite{kitamura2016determined,kitamura2018}.

In BSS/BSE, there is a wide class of methods that are based on Gaussian statistical models of signals, as compared to the non-Gaussianity-based methods considered in this paper. Those methods exploit only second order statistics (SOS) and their algebraic properties. For example, the analogy of the standard ICA problem based on SOS boils down to the problem of Joint Approximate Diagonalization (JAD) of covariance matrices; see, e.g., \cite{pham2001,tichavsky2009,yeredor2010} and the references therein. Similarly to IVA, the SOS-based methods were considered in \cite{li2009, anderson2012}; see also \cite{lahat2016}. 

\subsubsection{Locally Determined Models for Underdetermined BSS}
When the mixing model \eqref{ICA} involves more sources than observations ($n>d$), the extraction/separation and the mixing matrix identification problems are no more equivalent. Therefore, they are typically treated separately in two step procedures. For example, the estimation of $\A$ can be done by applying a decomposition to a tensor that is built from covariance matrices \cite{lathauwer2008} or higher-order based statistics \cite{ferreol2005,comon2006}. Then, various array processing methods can be applied to extract the sources \cite{tichavsky2011TSP,koldovsky2013SPL}.

There are also BSS methods that treat the underdetermined problem by assuming a certain local condition that guaranties that the every sample or time-frequency point involves maximally $d$ sources. Most typically, blind speech separation methods exploit the time-frequency sparsity of speech signals \cite{yilmaz2004,araki2004icakyoto}. Other methods assume that there are single-source points or regions and the separation mainly relies on a detection of these regions \cite{abrard2005,li2006}. Locally determined mixing is considered, e.g., in \cite{liu2014}. 

The CMV and CSV models, respectively, described through \eqref{CMV} and  \eqref{CSV} could be considered as members of the class of locally determined models for BSE, where the identification and extraction proceed jointly.

\subsubsection{Performance bounds}
Performance limitations of ICA based on the standard determined mixing model have been well investigated in the literature. It is known that $\A$ in \eqref{ICA} can be identified up to the order and scales of its 
columns if it holds that at most one source signal has the complex
Gaussian pdf or that no two complex Gaussian source signals have the same 
circularity coefficient \cite{complexErikson}. Then, a de-mixing matrix ${\bf W}$ can be estimated 
as such that ${\bf G}={\bf W}{\bf A}\approx{\bf P}\boldsymbol{\Lambda}$, where 
${\bf P}$ and $\boldsymbol{\Lambda}$ is a permutation and diagonal matrix (with nonzero diagonal entries), respectively.
${\bf G}$ reflects the separation accuracy as its $ij$th element, $G_{ij}$, determines the presence of $u_j$ in the $i$th separated signal $y_i$, so there is a clear correspondence between the elements of ${\bf G}$ and the Interference-to-Signal Ratio (ISR) of the separated signals. %The minimum possible value of ${\rm E}[|G_{ij}|^2]$ provides a lower bound for the accuracy achievable through the given model. 
For the real-valued (and similarly for the complex-valued) ICA problem, it was derived using the CRLB that the ISR of the $i$th separated source obeys %assuming that $u_1,\dots,u_d$ have unit variance (without any loss on generality), 
\begin{equation}\label{CRLB}
{\rm E}[\isr_i]\geq\frac{1}{N}\sum_{j=1, j\neq i}^d\frac{\okappa_j}{\okappa_i\okappa_j-1},
\end{equation}
where $N$ is the number of i.i.d. samples \cite{tichavsky2006,loeschCRB}; $\kappa_i={\rm E}[|\psi_i|^2]$ where $\psi_i(x)=-\partial/\partial x\,\log p_i(x)$ is the score function related to $p_i$, and $\okappa_i=\kappa_i\sigma^2_i$ where $\sigma^2_i$ is the variance of $u_i$; $\okappa_i$ corresponds to $\kappa_i$ when $p_i$ is normalized to unit variance. It holds that $\okappa_i\geq 1$, and $\okappa_i=1$ if and only if the $i$th pdf is circular Gaussian. Hence, the denominator in \eqref{CRLB} approaches zero when both the $i$th and the $j$th source signals are close to circular Gaussian.

This brings some things into question regarding the BSE problem. Without loss on generality, let $d-1$ source signals in the mixture be circular Gaussian but not so the first source (SOI). Then, $\A$ is no more identifiable, and the CRLB \eqref{CRLB} formally does not exist. However, BSE methods exploiting the non-Gaussianity of the SOI are known for their ability to blindly extract that source; see, e.g.,\cite{tichavsky2006}. Moreover, their asymptotic performance analyses have shown that their accuracy is limited by
\begin{equation}\label{CRLBICE}
{\rm E}[\isr]\geq\frac{1}{N}\frac{d-1}{\okappa-1}, \qquad i=2,\dots,d,
\end{equation}
where %${\bf q}^T={\bf w}^T\A$, ${\bf w}$ is the separation vector, and 
$\okappa=\okappa_1$; see, e.g., \cite{hyvarinen1997b,pham2006,tichavsky2006}. This asymptotic bound coincides with the right-hand side of \eqref{CRLB} when considering $i=1$ and $\okappa_j=1$ for $j=2,\dots,d$.

A formal confirmation of this bound for the real-valued case has been proven recently in \cite{kautsky2017} through computing the CRLB for the ICE mixing model, that is, assuming that the mixing matrix is structured as described by \eqref{mixing} and that the background signals are Gaussian. 

In the first part of this paper, we generalize this result for the complex-valued case where the SOI is assumed to be non-Gaussian or non-circular Gaussian. The background is modeled as circular Gaussian or circular non-Gaussian. We avoid the case with non-circular background, for simplicity, as it is computationally less tractable and its analysis goes beyond the scope of this paper. We show that the CRLB of ICE corresponds with the bound for ICA when the background is circular Gaussian, as in the real-valued case. Moreover, we also show that these bounds coincide when the background modeling in ICE takes into account possible non-Gaussianity of the background. %Moreover, it can be even higher than that of ICA when the background latent sources are dependent such that they cannot be linearly decomposed into independent sources.

%\subsection{Contribution}
%BSS has been solved, in the context of ICA, through extracting a signal by minimizing its entropy \cite{comon2010handbook}. The information theory-based approach, however, does not allow us to directly compute the CRLB for the extraction problem. Recently, we have revised the BSE in approach called Independent Component Extraction (ICE) \cite{eusipco2017}. Here, the mixing model \eqref{ICA} is re-parameterized so that it contains a minimum number of parameters that are necessary for the extraction of the SOI. The statistical model is based on the assumption that the SOI is independent from the background, similarly to ICA, and the background is assumed to be Gaussian. Then, the computation of the CRLB is straightforward using the likelihood function. 

%This way, we have computed the bound for the real-valued case in \cite{kautsky2017}. It was shown that the bound coincides with that for ICA if the background is Gaussian and is also in a good agreement with asymptotic performance analyses of several BSE algorithms \cite{tichavsky2006,hyvarinen1997b}. In this paper, we generalize this result for the complex-valued case where the SOI is assumed to be non-Gaussian or non-circular Gaussian, while the background is modeled as circular Gaussian. Analysis of complex-valued case is an important extension, since 

The article is organized as follows. Section~II is devoted to the standard determined mixing model and the above mentioned issues related to the CRLBs. In Section~III the piecewise determined mixing models are introduced, and the related CRLBs are derived using results of Section~II. The computed theoretical bounds are discussed and compared in Section~IV through analyzing several special cases. Experimental validations of the bounds are presented in Section~V, and the article is concluded by Section~VI.

\subsection{Nomenclature} 
Plain letters denote scalars, bold letters denote 
vectors, and bold capital letters denote matrices. Upper index $\cdot^T$, 
$\cdot^H$, or 
$\cdot^*$ denotes, respectively, transposition, conjugate transpose, or complex 
conjugate. The Matlab convention for matrix/vector concatenation and 
indexing will be used, e.g., $[1;\,{\bf g}]=[1,\, {\bf g}^T]^T$, and 
$(\A)_{j,:}$ is the $j$th row of $\A$. 
A complex random vector $\x$ is called circular if its pseudo-covariance is
$\mathtt{pcov}(\x)={\rm E}\bigl[\left(\x-{\rm E}[\x]\right)\left(\x-{\rm 
E}[\x]\right)^T\bigr]={\bf 0}$, otherwise, $\x$ is non-circular; ${\rm 
E}[\cdot]$ stands for the expectation operator. The second-order 
circularity coefficient $\gamma$ of a complex-valued random variable $x$ with zero mean, see \cite{complexErikson}, is 
defined as in \cite{cGGD} $\rho = \left|{\rm E}\bigl[x^2\bigr]\right|/{\rm E}\bigl[\left|x\right|^2\bigr]$. Thus, $\rho \in [0,1]$ and $\rho = 0 $ holds for circular random variable.

\section{Determined Mixing}\label{secmodel}
\subsection{Algebraic Model}
%Without any loss of generality, let the SOI be $s=u_1$. 
Here, we briefly explain the parameterization of \eqref{ICA} as given by \eqref{mixing} and \eqref{demixing}. 
Let $\A$ be partitioned as $\A=[{\bf a},\,\A_2]$. Then, 
$\x$ can be written as 
\begin{equation}\label{definition_y}
 \x=\A{\bf u}=\a s+{\bf y}, 
\end{equation}
where ${\bf y}=\A_2{\bf u}_2$ 
and ${\bf u}_2=[u_2,\dots,u_d]^T$. Since neither ${\bf u}_2$ nor $\A_2$ should be estimated in order to extract $s$, we can consider any auxiliary background signals $\z$ such that ${\bf y}=\A_2{\bf u}_2 = {\bf Q}\z$, where the columns of ${\bf Q}$ span the same subspace as those of $\A_2$. Compared to ${\bf u}_2$, the elements of $\z$ need not be independent, so ${\bf Q}$ can be arbitrary in this sense. 

The structures \eqref{mixing} and \eqref{demixing} are obtained based on the following three conditions
\begin{align}
{\bf B}{\bf a} &=  {\bf 0},\\
{\bf w}^H{\bf Q} &= {\bf 0}^T,\\
\W_{\rm ICE}\A_{\rm ICE} &= {\bf I}_d,
\end{align}
where the first two conditions are, in fact, involved in the third one.
These conditions ensure that ${\bf w}^H \x = s$ and ${\bf B}\x = \z$, in other words, that $\W_{\rm ICE}$ is de-mixing, i.e., it extracts $s$ from $\x$ and separates it from $\z$. The ICE algebraic model can thus be written as 
\begin{equation}\label{mixingICE}
\x = \A_{\rm ICE}{\bf v},
\end{equation}
where ${\bf v}=[s;{\bf z}]$. 

\subsection{Statistical Model}

%Several stochastic models have been considered in BSS/ICA to model the 
%signals' independence, e.g., relying on one or more signal properties such as 
%non-Gaussianity, noncircularity, nonstationarity or nonwhiteness 
%\cite{adali2014}. In ICE, there are only two complex variables: $s$, which 
%coincides with $u_1$, and ${\bf z}$, which is a vector variable having an 
%unspecified structure (it is a mixture of $u_2,\dots,u_d$).

The fundamental assumption of ICA/ICE states that $s$ and $\z$ are
independent, which means that their joint pdf can be factorized as the product 
of marginal pdfs. Let the pdfs of $s$ and $\z$ be denoted $p_s(s)$ and $p_\z(\z)$, respectively.
% while $\z$ has multivariate circular Gaussian pdf with 
%covariance $\C_{\bf z}$. 
%The latter assumption can be justified by the fact that, up to very special 
%cases, $\z$ is a mixture of ${\bf u}_2$. Even if $u_2,\dots,u_d$ are 
%non-Gaussian circular (eventually noncircular), their mixture tends to have 
%distribution close to circular (eventually noncircular) Gaussian due to the 
%Central Limit Theorem \cite{adaliComplexSP}. %The ICA and the ICE models 
%%coincide when $u_2,\dots,u_d$ are circular Gaussian.
Using \eqref{mixingICE}, the pdf of $\x$ is
\begin{equation}\label{p_x}
   p_{\x}(\x)=p_{s}({\bf w}^H\x)p_\z({\bf B}\x)|\det(\W_{\rm ICE})|^2,
\end{equation} 
where $\det(\W_{\rm ICE})=(-1)^{d-1}\gamma^{d-2}$. 
%As shown later, the de-mixing matrix can be scaled to have $|\det(\W_{\rm ICE})|=1$. 

\subsection{Indeterminacies}\label{sec:indeterminacies}
ICE involves that same indeterminacies as ICA as the problem is solved through finding vector parameters ${\bf w}$ and ${\bf a}$ such that $s$ and $\z$ are independent. It follows that 
any independent component of $\x$ could play the role of $s$, because of the indeterminacy of the order of original components in \eqref{ICA}. In this work, this problem can be overlooked as the CRLB analysis is local. In practice, any estimating algorithm must be properly initialized in order to extract the desired source.

The scales of $s$ and of ${\bf a}$ are ambiguous in the sense that $s$ and ${\bf a}$ can be substituted, 
respectively, by $\alpha s$ and $\alpha^{-1}{\bf a}$ with any $\alpha\neq 0$. This is know as the scaling ambiguity problem. Since Interference-to-Signal Ratio is invariant to the scaling, we can later cope with this ambiguity by fixing some scalar parameter in the mixing model. In this section, we put $\gamma=1$.

% Then, the only free 
%parameters of the ICE model are ${\bf g}$ and ${\bf h}$.

%where $\W_{\rm ICE}$, ${\bf w}$, and ${\bf B}$ depend on ${\bf g}$ and ${\bf 
%h}$ as described by \eqref{demixing}, and $p_\z$ is the pdf of the circular 
%Gaussian distribution with zero mean and covariance $\C_{\bf z}$, denoted by 
%$\mathcal{CN}({\bf 0},\C_\z)$. 

%A straightforward calculus, not shown here to save space, can show that for 
%$\gamma=1$, $|\det(\W_{\rm ICE})| = 1$; see \cite{eusipco2017}. The 
%log-likelihood function for one signal sample is thus equal to
%\begin{multline}\label{likelihood}
%\mathcal{L}(\x|\boldsymbol{\theta},\boldsymbol{\xi} ) = \log p_{s}({\bf w}^H\x)-\x^H {\bf 
%B}^H {\textbf{C}_{\z}^{-1}} {\bf B}\x - \log\left( \det(
%\textbf{C}_{\z})\right) \\ - (d-1)\log(2\pi),
%\end{multline} 
%where $\boldsymbol{\theta} = [\g;\h;\bf r]$ denotes the parameter vector, in which 
%$\bf r$ stands for the $d(d-1)/2\times 1$ vector stacking the off-diagonal 
%elements of $\textbf{C}_z^{-1}$, and $\boldsymbol{\xi}$ is the $d\times 1$ 
%vector stacking the real-valued diagonal elements of $\textbf{C}_z^{-1}$; 
%$|\textbf{C}_{\z}|$ denotes the determinant of 
%$\textbf{C}_{\z}$; ${\bf r}$ and $\boldsymbol{\xi} $ are nuisance parameters. 

\subsection{Interference-to-Signal Ratio}
%A lower bound for the achievable mean Interference-to-Signal Ratio (ISR) can be derived using the CRLB. 
Let $\widehat{\bf w}$ be an estimated separating vector ${\bf w}$. %and let ${\bf A}_{\rm ICE}$ be the true mixing matrix in \eqref{mixingICE}. 
Using \eqref{definition_y}, the extracted signal is equal to $\widehat{s}=\widehat{\bf w}^H\x=\widehat{\bf w}^H{\bf a} + \widehat{\bf w}^H{\bf y}=\widehat{\bf w}^H{\bf a} + \widehat{\bf w}^H{\bf Qz}$. The ISR of the signal is 
\begin{equation}\label{ISRinst}
\isr=\frac{{{\rm E}}[|\widehat{\bf 
w}^H\y|^2]}{{{\rm E}}[|\widehat{\bf w}^H{\bf a}s|^2]}=
\frac{{\bf q}_2^H\C_\z{\bf q}_2}{|q_1|^2\sigma_s^2} \approx \frac{1}{\sigma_s^2}{\bf q}_2^H\C_\z{\bf q}_2,
\end{equation} 		
where ${\bf q}^H=[q_1,\,{\bf q}_2^H]= %\widehat{\bf w}^H{\bf A}_{\rm ICE}=
[\widehat{\bf w}^H{\bf a},\,\widehat{\bf w}^H{\bf Q}]$, and $\C_\z$ stands for the covariance matrix of $\z$. The last approximation in \eqref{ISRinst} is valid for ``small'' estimation error in $\widehat{\bf w}$, that is, when ${\bf q}\approx{\bf e}_1$ (the unit vector). Then, the mean ISR value reads  
\begin{equation}\label{ISRmean}
 {{\rm E}}\left[\isr\right]\approx 
 \frac{1}{\sigma_s^2}{{\rm E}}\left[{\bf q}_2^H\C_\z{\bf q}_2\right] = 
 \frac{1}{\sigma_s^2}{\tt tr}\left(\C_\z{{\rm E}}\left[{\bf q}_2{\bf 
 q}_2^H\right]\right).
\end{equation}  
Hence, \eqref{ISRmean} can be written as
\begin{equation}\label{ISRapro}
 {\rm E}\left[\isr\right]\approx \frac{1}{\sigma_s^2}{\tt 
 tr}\left(\C_\z\cov\left({\bf q}_2\right)\right),
\end{equation} 
where we can see that the covariance matrix of ${\bf q}_2$, denoted as $\cov\left({\bf q}_2\right)$, characterizes the accuracy of $\widehat{\bf w}$. By replacing $\cov\left({\bf q}_2\right)$ by the corresponding CRLB, we obtain the algorithm-independent Cram\'er-Rao-induced bound (CRIB) for ISR.

%To compute the covariance matrix of ${\bf q}_2$, we perform a parametrization of the mixing matrix $\A_{\rm ICE}$ in \eqref{mixingICE}, see \cite{eusipco2017}. Without any loss of generality, let the first column of $\A_{\rm ICE}$ be $\a = [1; \g]$. Then it holds, $\y = \A_2 {\bf u}_2 = {\bf Q}\z$ and obviously ${\bf B} \A_2 {\bf u}_2 = {\bf B}{\bf Q}\z = \z$ holds for any ${\bf Q}$, for which exists a matrix ${\bf B}$ such as ${\bf B}{\bf Q}$ is non-singular (there is a bijection between the original signal ${\bf u}_2$ and $\z$). Useful selection, discussed in \cite{eusipco2017}, is 
%\begin{equation}\label{Qparam}
%{\bf Q} = \begin{pmatrix}
%\h^H \\
%\g\h^H-{\bf I}_{d-1}
%\end{pmatrix}.
%\end{equation} Using \eqref{Qparam} enables the computation of ${\bf B}$ and ${\bf w}$. Straitforward calculation shows 
%\begin{equation}
%{\bf B} = \left[\g, \ -{\bf I}_{d-1}\right], \quad {\bf w} = \begin{pmatrix}
%1-\g^H\h \\
%\h
%\end{pmatrix}.
%\end{equation} 
%The equivariance property of the BSE problem \cite{tichavsky2006} enables us to consider the special case when $\g=\h=\0$, without any loss on generality.

\subsection{Cram\'er-Rao-induced Bound}\label{section:CRIB}
%To derive the bound, the Fisher information matrix of ${\bf q}_2$ will be computed. 
Let the parameter vector be $\boldsymbol{\theta} = [{\bf a}; {\bf w}]$. In the following, we exploit a transformation rule saying that the Fisher Information Matrix (FIM) of $\boldsymbol{\theta}$, denoted as ${\bf F}_{\boldsymbol{\theta}}$, and the FIM of a linearly transformed version $\boldsymbol{\varphi} = {\bf K}\boldsymbol{\theta}$, where ${\bf K}$ is a non-singular matrix, are related through \cite{complCRB}
\begin{equation}\label{Fphi}
     \textbf{F}_{\boldsymbol{\varphi}}= \textbf{K}^{-1}\textbf{F}_{\boldsymbol{\theta}}\textbf{K}^{-H}.
\end{equation}
This property will be used to show that we can derive the CRIB for \eqref{ISRapro} by considering CRLB when the mixing parameters are $\h={\bf 0}$. This property is related to the {\em equivariance} of the BSS mixing model \eqref{ICA}, see, e.g., \cite{laheld1996,comon2010handbook}.

Now, consider the special case when ${\bf h} = {\bf g}= {\bf 0}$, for which the parameter vector is equal to $\boldsymbol{\theta}_{\bf I} = [{\bf e}_1; {\bf e}_1]$. The transform between $\boldsymbol{\theta}$ and $\boldsymbol{\theta}_{\bf I}$ is given by 
\begin{equation}\label{thetaI}
    \boldsymbol{\theta} = 
\underbrace{\begin{pmatrix}
	{\bf A}_{\rm ICE} & {\bf 0}\\
	{\bf 0} & {\bf W}_{\rm ICE}^H
\end{pmatrix}}_{\bf K} \boldsymbol{\theta}_{\bf I} = {\bf K}\boldsymbol{\theta}_{\bf I},
\end{equation}
where ${\bf A}_{\rm ICE}$ and ${\bf W}_{\rm ICE}$ are, respectively, given by \eqref{mixing} and \eqref{demixing}. According to \eqref{Fphi}, it holds that 
\begin{equation}\label{FEye}
\textbf{F}_{{\boldsymbol{\theta}}} = \textbf{K} \textbf{F}_{\boldsymbol{\theta}_{\bf I}} \textbf{K}^{H}.
\end{equation}
Similarly, we can consider a transformed parameter vector
\begin{equation}\label{phiEye}
    \boldsymbol{\theta}_{\bf q} = \underbrace{\begin{pmatrix}
{\bf W}_{\rm ICE} & {\bf 0}\\
{\bf 0} & {\bf A}_{\rm ICE}^H
\end{pmatrix}}_{{\bf K}^{-1}}
\begin{pmatrix}
 {\bf a}\\
{\bf w}
\end{pmatrix} = {\bf K}^{-1}\boldsymbol{\theta},
%=
%\begin{pmatrix}
% {\bf a}\\
%{\bf q}
%\end{pmatrix}.
\end{equation} 
and holds that $\textbf{F}_{{\boldsymbol{\theta}}_{\bf q}}= \textbf{K}^{-1}\textbf{F}_{\boldsymbol{\theta}}\textbf{K}^{-H}$, which, together with \eqref{FEye}, results in
\begin{equation}\label{equivariance}
    \textbf{F}_{{\boldsymbol{\theta}}_{\bf q}} = \textbf{F}_{\boldsymbol{\theta}_{\bf I}}.
\end{equation}
%Next, if $\varphi = \varphi(\theta)$ is a differentiable function of $\theta$, then the
%Fisher information matrix for exists as well and is equal to
%\begin{equation}\label{FIM transformation}
%    \textbf{F}_{\varphi}= \textbf{J}^{-1}\textbf{F}_{\theta}\textbf{J}^{-H}
%\end{equation}, where $\textbf{J}$ is the Jacobian of the mapping $\varphi(\theta)$. 
%In \eqref{phiEye} the mapping is $\varphi_{\bf q}(\boldsymbol{\theta}) = {\bf M}\boldsymbol{\theta}$. Thus by using \eqref{FIM transformation} \begin{equation}\label{Fphi}
%     \textbf{F}_{\varphi}= \textbf{M}^{-1}\textbf{F}_{\theta}\textbf{M}^{-H}.
%\end{equation}

%For ${\bf h} = {\bf 0}$, it holds that ${\bf q}_2=-\widehat{\bf h}$ where $\widehat{\bf h}=\widehat{\bf w}_{2:d}$.
Hence, from \eqref{equivariance} it follows that the CRIB for \eqref{ISRapro} can be obtained by replacing $\cov\left({\bf q}_2\right)$ by the corresponding CRLB, which is equal to the CRLB for the unbiased estimation of ${\bf h}$ when its true value is ${\bf h} = {\bf 0}$. Finally,
%Then, ${\bf q}_2=\widehat{\bf h}$, where $\widehat{\bf h}=\widehat{\bf w}_{2:d}$, and \eqref{ISRapro} simplifies to 
\begin{equation}\label{ISRvsCRLB}
 {{\rm E}}\left[\isr\right] \approx \frac{1}{\sigma_s^2}{\tt 
 tr}(\C_\z\cov(\widehat{\bf h})) \geq \frac{1}{\sigma_s^2}{\tt 
 tr}\left(\C_\z\crb\left(\h\right)|_{\h={\bf 0}}\right), 
\end{equation} 
where $\crb(\h)|_{\h={\bf 0}}$ denotes the diagonal block of the inverse matrix of the FIM 
%for ${\boldsymbol{\theta}}_{\bf I}$.
%$\mathcal{J}^{-1}\left(\boldsymbol{\tilde{\theta}}\right)$ 
corresponding to the parameter vector $\h$ when ${\bf h} = {\bf 0}$. 
%Hence, to complete the calculation of CRIB, the inverse of Fisher Information Matrix is needed. 

\subsection{Fisher information matrix} 
To compute the CRLB, we use the approach for the complex-valued parameters described in \cite{complCRB}. %Let 
%$\boldsymbol{\theta}$ and $\boldsymbol{\xi}$ denote, respectively, the complex 
%parameter vector and the real parameter vector to be estimated. 
By putting $\gamma=1$, as justified in Section~\ref{sec:indeterminacies}, the only free parameters of the mixing model \eqref{mixingICE} are $\h$ and $\g$, so let the parameter vector be $\boldsymbol{\theta} = 
[\h;\g]$. According to \cite{complCRB}, for any unbiased estimator of $\boldsymbol{\theta}$, 
it holds that
\begin{equation}
\cov\left( \boldsymbol{\theta}\right) \succeq 
\mathcal{J}^{-1}\left(\boldsymbol{\theta} \right) = 
\crb\left(\boldsymbol{\theta} \right),
\end{equation} where $\mathcal{J}\left(\boldsymbol{\theta} \right)$ 
is the FIM, and ${\bf C}\succeq {\bf D}$ means that 
${\bf C}-{\bf D}$ is a positive semi-definite matrix. $\mathcal{J}\left(\boldsymbol{\theta} \right)$ can be 
partitioned as 
\begin{equation}\label{complexFIM}
\mathcal{J}\left(\boldsymbol{\theta} \right) = \begin{pmatrix}
 \textbf{F} & \textbf{P}  \\
 \textbf{P}^* &  \textbf{F}^*  
  \end{pmatrix},
\end{equation} 
where % ${\bf F}_{\boldsymbol{\theta}}$, ${\bf P}_{\boldsymbol{\theta}}$, $ 
%\textbf{R}_{\boldsymbol{\theta}}$ and $\textbf{F}_{\boldsymbol{\xi}}$ are 
%defined as 
  \begin{equation}\label{wirtingerDerivatives}
 \textbf{F} =
 {{\rm E}}\left[\frac{\p \mathcal{L}}{\p \boldsymbol{\theta}^*}\left(\frac{\p 
 \mathcal{L}}{\p 
 \boldsymbol{\theta}^*}\right)^H\right], \quad
  \textbf{P}  = 
 {{\rm E}}\left[\frac{\p \mathcal{L}}{\p \boldsymbol{\theta}^*}\left(\frac{\p 
 \mathcal{L}}{\p 
 \boldsymbol{\theta}^*}\right)^T\right] 
%\label{RFIM}
% \textbf{R} & = 
% {{\rm E}}\left[\frac{\p \mathcal{L}}{\p \boldsymbol{\theta}^*}\left(\frac{\p 
% \mathcal{L}}{\p 
% \boldsymbol{\xi}}\right)^T\right], &
% \textbf{F}_{\boldsymbol{\xi}} & = 
% {{\rm E}}\left[\frac{\p \mathcal{L}}{\p \boldsymbol{\xi}}\left(\frac{\p 
% \mathcal{L}}{\p 
% \boldsymbol{\xi}}\right)^T\right],
  \end{equation}
and where the derivatives in \eqref{wirtingerDerivatives} are defined 
according to the Wirtinger calculus. $\mathcal{L}(\cdot)$ denotes the log-likelihood function of \eqref{p_x}, namely, 
\begin{equation}\label{likelihood}
 \mathcal{L}(\h,\g|\x) = \log p_{s}({\bf w}^H\x)+ \log p_\z({\bf B}\x)
\end{equation}
%  An additional assumption follows the equivariance property \cite{tichavsky2006}: we can assume the true values of unknown parameters equal to zero, i.e. $\h=\g={\bf 0}$. Let the symbol $|_{\h=0}$ denotes an evaluation in $\h=\vec{0}$.
 
The derivatives of the log-likelihood function \eqref{likelihood} are as 
follows:
\begin{align} \label{dg}
\frac{\p \mathcal{L}(\x|\boldsymbol{\theta})}{\p \g^*}\Big|_{\h=0} &= -\boldsymbol{\psi}_{\z}(\z)s^* , \\ \label{dh}
\frac{\p \mathcal{L}(\x|\boldsymbol{\theta})}{\p \h^*}\Big|_{\h=0} &= \psi_s^*(s)\z ,
 \end{align} 
 where  $\psi_s(s) = 
 -\frac{\p \ln p_s(s,s^*)}{\p s^*}$ and $\boldsymbol{\psi}_{\z}(\z) = 
 -\frac{\p \ln p_{\z}(\z,\z^*)}{\p \z^*}$ are the score functions. Using \eqref{dg},\eqref{dh}, $\textbf{F}$ in \eqref{wirtingerDerivatives} is calculated as
 \begin{equation}\label{realFim} 
 \textbf{F}=%(\boldsymbol{\theta})=  \textbf{F}(\g,\h,\bf r)= 
 \begin{pmatrix}
 \sigma_s^2 \boldsymbol{\kappa}_{\z} & -\textbf{I}_{d-1}  \\
 -\textbf{I}_{d-1} &  \kappa_s  \textbf{C}_{\z}
  \end{pmatrix}. 
\end{equation} 
 where \begin{align}
 \kappa_s &= {\rm E}[|\psi(s)|^2],\\
 \sigma_s^2 &= {\rm E}[|s|^2], \\
 \bkappa_{\z} &= {\rm E}\left[\boldsymbol{\psi}_{\z}(\z) \boldsymbol{\psi}_{\z}^H(\z) \right]\label{kappaz}
  \end{align}   

  Now, we describe the computation of $\textbf{P}$ in \eqref{complexFIM}. Let 
  $\textbf{P}$  be partitioned as
 \begin{equation}
 \textbf{P}=%(\boldsymbol{\theta})=  \textbf{P}(\g,\h,\bf r)= 
 %\left(\begin{array}{ccc}
 \begin{pmatrix}
 \textbf{P}_{\g,\g} & \textbf{P}_{\g,\h}  \\
 \textbf{P}_{\g,\h}^T &  \textbf{P}_{\h,\h}
  \end{pmatrix}. %\right). 
\end{equation} 
Then,
 \begin{align}
 \textbf{P}_{\g,\g} =&\  {\rm E}\left[\boldsymbol{\psi}_{\z}(\z)\boldsymbol{\psi}_{\z}^T(\z)\right]{\rm E}\left[{s^*}^2 \right],  \\
   \textbf{P}_{\h,\h} =&\  {\rm 
  E}\left[\psi_s^*(s)^2\right] {\rm E}\left[\z\z^T\right], \\
  \textbf{P}_{\g,\h} =&\ {\bf 0}.
 % {\rm E}\left[s^*\psi_s^*(s)\right]{\rm E}\left[\boldsymbol{\psi}_{\z}(\z)\z^T\right].
\end{align} 
%Calculation of CRIB? 

\subsection{Circular sources}\label{sec:circularsources}
In general, the analytic computation of the inverse matrix of \eqref{complexFIM} is not tractable. Therefore, we investigate two special cases in the following subsections. 

Here, we assume that $s$ and $\z$ have general {\em circular} pdf. Assuming this, the FIM \eqref{complexFIM} obtains the block-diagonal structure, because $\textbf{P}_{\h,\h} = \textbf{P}_{\g,\h} = {\bf 0}$ due to the circularity of $\z$ and $\textbf{P}_{\g,\g} = {\bf 0}$ due to the circularity of $s$, and, then,
  \begin{equation}\label{FIMproCirculars}
\mathcal{J}\left(\boldsymbol{\theta}\right)=\left(\begin{array}{ccc}
 \sigma_s^2 \boldsymbol{\kappa}_\z^{-1} &  -\I_{d-1} & \bf {O} \\
 -\I_{d-1} &  \kappa_s \C_\z & \bf {O} \\
 \bf {O} & \bf {O} &  {\textbf{F}^*} \\ 
  \end{array}\right).
\end{equation} 
%where ${\bf D}$ is not the target of interest to derive in this paper.
$\crb(\h)|_{\h={\bf 0}}$ is obtained as the upper right diagonal block of the inverse matrix of \eqref{FIMproCirculars}, which reads  
 \begin{equation}\label{CRBcirculars}
\crb(\h)|_{\h={\bf 0}} = \left({\kappa_s}\C_\z-\frac{1}{\sigma_s^2}\boldsymbol{\kappa}_\z^{-1}\right)^{-1} .
\end{equation} 
Applying the transformation theorem in \eqref{kappaz}, it can be shown that, for $\tilde\z = {\bf T}\z$, it holds that
\begin{equation}\label{transfkappa}
    \bkappa_{\bf z} = {\bf T}{\bkappa}_{\tilde\z}{\bf T}^H,
\end{equation} 
where ${\bf T}$ is a non-singular transformation matrix. By taking ${\bf T}={\bf C}_{\bf z}^{-\frac{1}{2}}$, which is a matrix satisfying that ${\bf C}_{\bf z}^{-\frac{1}{2}}{\bf C}_{\bf z}^{-\frac{1}{2}}={\bf C}_{\bf z}^{-1}$, then $\bkappa_{\tilde\z}$ corresponds to the statistic of uncorrelated and unit-scaled $\z$. Hence, \eqref{CRBcirculars} can be written as
\begin{equation}\label{CRBcirculars2}
\crb(\h)|_{\h={\bf 0}} = \C_\z^{-\frac{1}{2}} \left({\kappa}_s{\bf I}_{d-1}-\frac{1}{\sigma_s^2}\boldsymbol{\bkappa}_{\tilde\z}^{-1}\right)^{-1} \C_\z^{-\frac{1}{2}}.
\end{equation}

By putting \eqref{CRBcirculars2} into \eqref{ISRvsCRLB}, %and using the permutation property of the matrix product trace, 
the CRIB for ISR, when considering $N$ observations, is
\begin{equation}\label{CRIBcirculars}
{\rm E}\left[\isr\right] \geq  \frac{1}{N} 
\frac{1}{\sigma_s^2}\tr \left[ \left({\kappa_s}{\bf I}_{d-1}-\frac{1}{\sigma_s^2}\boldsymbol{\kappa}_{\tilde{\z}}^{-1}\right)^{-1}\right].
\end{equation} 

Next, we can use the identity \eqref{transfkappa} again by considering ${\bf T}$ such that ${\bf T}\tilde\z$ are independent. Since $\tilde\z$ are uncorrelated and normalized, such ${\bf T}$ must be unitary, i.e., ${\bf T}{\bf T}^H={\bf I}_{d-1}$. Also, provided that all but one components in the original model \eqref{ICA} are non-Gaussian, the entire mixture is separable, so ${\bf T}\tilde\z$ must be equal to ${\bf u}_2$ up to their order and scales. Without any loss on generality, we can assume that ${\bf T}$ is such that ${\bf T}\tilde\z={\bf u}_2$ and that ${\bf u}_2$ have unit variance.  Then, ${\bkappa}_{{\bf T}\tilde\z}$ is diagonal having %\begin{equation}
%    \tilde{\kappa}_\z = {\bf T}^{-1}{\bf D}{\bf T}, 
%\end{equation}
%where ${\bf D}$ is a diagonal matrix with eigenvalues of $\tilde{\kappa}_\z$ on its diagonal, denoted as 
diagonal elements equal to $\okappa_2,\dots,\okappa_{d}$, and \eqref{CRIBcirculars} simplifies to
\begin{equation}\label{cribEigenVals}
    {\rm E}\left[\isr\right] \geq  \frac{1}{N} 
\sum_{j=2}^{d} \frac{\okappa_j}{\sigma_s^2{\kappa_s}{\okappa_j}-1}.
\end{equation}
This bound corresponds with \eqref{CRLB} for $i=1$ since $\sigma_s^2\kappa_s=\okappa_s=\okappa_1$, which means that the same extraction accuracy can be achieved through ICE as that by ICA. It should be, however, noted that the multivariate score function $\psi_\z(\cdot)$ must be known for realizing maximum likelihood estimation \cite{koldovsky2018a}. 

In our considerations, we can go also slightly beyond the standard ICA. Let the observed signals obey the model \eqref{mixingICE} but not \eqref{ICA}, that is, there need not exist ${\bf T}$ such that ${\bf T}\tilde\z$ are independent (no independent components $u_2,\dots,u_d$ are assumed, only the independence between $s$ and $\z$). Since ${\bkappa}_{\tilde\z}$ is positive definite, we can consider its decomposition
\begin{equation}
    \bkappa_{\tilde\z}={\bf U}{\bf D}{\bf U}^H
\end{equation}
where ${\bf U}^H$ is the unitary matrix of eigenvectors of ${\bkappa}_{\tilde\z}$, and ${\bf D}$ is diagonal with diagonal entries denoted as $\omega_2,\dots,\omega_d$.
%It holds that 
%\begin{align}
%    \tr(\tilde{\kappa}_\z) = &\sum_{i=1}^{d-1}\lambda_i = \sum_{i=1}^{d-1}\omega_i\\
%    \sum_{i=1}^{d-1}\frac{1}{\omega_i} \leq& \sum_{i=1}^{d-1}\frac{1}{\kappa_i}. \label{eigenValInequality}
%\end{align}
%Then, \eqref{cribEigenVals} could be written as 
%\begin{equation}
%    \frac{1}{N} \sum_{i=1}^{d-1} \frac{\lambda_i}{\sigma_s^2{\kappa_s}{\lambda_i}-1} = %\frac{1}{N}\frac{1}{\sigma_s^2{\kappa_s}}\sum_{i=1}^{d-1}
%    \left(1 + \frac{1}{\sigma_s^2{\kappa_s}{\lambda_i}-1} \right).
%\end{equation} 
%Thus, combining \eqref{cribEigenVals} and \eqref{eigenValInequality} leads to 
Then, the CRIB obtains similar form to \eqref{cribEigenVals}
\begin{equation}\label{CRLB_ngICE}
    {\rm E}\left[\isr\right]  \geq  \frac{1}{N} 
\sum_{j=2}^{d} \frac{\omega_j}{\sigma_s^2{\kappa_s}{\omega_j}-1}. %\geq  \frac{1}{N} 
%\sum_{i=1}^{d-1} \frac{\lambda_i}{\sigma_s^2{\kappa_s}{\lambda_i}-1}, 
\end{equation} 

%which shows that the CRIB for ICE could be even lower then that one for ICA thanks to the dependencies between the background signals. 

\subsection{Circular Gaussian Background}\label{section:gausscase}
Here, we assume that $s$ can be arbitrary non-circular and non-Gaussian while $\z$ is circular Gaussian. Under this assumption, $\textbf{P}_{\h,\h} = {\bf 0}$, and
%Since $\z = {\bf B} \y = {\bf B}{\bf A}_2{\bf u}_2 $, then
%if the original signals, ${\bf u}_2$, are circular, then $\z$ is circular.  Moreover, due to the Central Limit Theorem (CLT) \cite{adaliComplexSP} the distribution of $\z$ tends to multivariate Gaussian for large number of sources. Considering the additional assumption of circularity and using the CLT simplify the FIM, 
since $\boldsymbol{\kappa}_{\z} = \C_{\z}^{-1}$, also $\textbf{P}_{\g,\g}= {\bf 0}$ thanks to the circularity of $\z$. The FIM thus obtains a similar structure to \eqref{FIMproCirculars}, namely,
\begin{equation}\label{FIMproI=Cz}
\mathcal{J}\left(\boldsymbol{\theta}\right)=\left(\begin{array}{ccc}
 \sigma_s^2 \C_\z^{-1} &  -\I_{d-1} & \bf {O} \\
 -\I_{d-1} &  \kappa_s \C_\z & \bf {O} \\
 \bf {O} & \bf {O} &  {\textbf{F}^*} \\ 
  \end{array}\right).
\end{equation} 
%where ${\bf D}$ is not the target of interest to derive in this paper.
%Since $\mathcal{J}\left(\boldsymbol{\theta}\right)$  is a block diagonal matrix, the block for $\crb(\h)$ reads  
Hence,
 \begin{equation}\label{CRB}
\crb(\h)|_{\h={\bf 0}} = \left({\kappa_s}\C_\z-\frac{1}{\sigma_s^2}\C_\z\right)^{-1} = \frac{\sigma_s^2}{\kappa_s\sigma_s^2-1}\C_\z^{-1}.
\end{equation} 
Therefore, for $N$ observations, the CRIB for ISR says that
\begin{equation}\label{CRIB}
{\rm E}\left[\isr\right] \geq  \frac{1}{N} 
\frac{d-1}{\okappa_s-1}.
%\tr \left(\textbf{I}_{d-1}\right) = \frac{d-1}{N} 
%%%\frac{1}{\kappa\sigma_s^2-1}\frac{\sigma_s^2}{\sigma_{\z}^2} ,
\end{equation}
This result confirms the asymptotic bound given by \eqref{CRLBICE} for complex-valued non-circular SOI. 

%, that the considered parametrization is not restrictive, since the result is in a good agreement with \eqref{CRLB} in the following sense: 
%Assume that $u_1,\dots,u_d$ have all unit variance and that $u_2,\dots,u_d$ 
%have all circular Gaussian pdf, which means that $\kappa_j=1$ for 
%$j=2,\dots,d$. In that special case, the bound for the achievable ISR through 
%ICA induced by \eqref{CRLB} coincides with \eqref{CRIB}. Since, $\kappa 
%\sigma_s^2 \approx \kappa_{norm}$, where $\kappa_{norm} = {\rm E}[|\psi(s)|^2]$ 
%when $s$ is normalized to the unit variance.

\section{Piecewise Determined Mixing}
We now turn to the piecewise determined mixtures, in general, described by \eqref{newmodel}.
For simplicity, let the number of available samples in each of $M$ blocks be the same, equal to $N_b$. It holds that $M\cdot N_b=N$. The variance of the SOI and the covariance matrix of the background signals in the $m$th block will be denoted by $\sigma_{s^m}^2$ and $\C_{\z^m}$, respectively.

Let $\widehat{\bf w}^m$ be an estimated separating vector for the $m$th block, $m=1,\dots,M$. The ISR of the extracted signal evaluated over the entire data is equal to
\begin{multline}\label{ISRblocks}
  \isr =  \frac{{\sum_{m=1}^M {\rm E}}[|(\widehat{\bf 
w}^m)^H\y^m|^2]}{\sum_{m=1}^M{{\rm E}}[|(\widehat{\bf w}^m)^H{\bf a}^ms^m|^2]}
= \frac{\sum_{m=1}^M({\bf q}_2^m)^H\C_{\z^m}{\bf q}_2^m}{\sum_{m=1}^M|q_1^m|^2\sigma_{s^m}^2} = \\
= \frac{\sum_{m=1}^M{\tt tr}\Bigl(\C_{\z^m}{\bf q}_2^m({\bf q}_2^m)^H\Bigr)}{\sum_{m=1}^M|q_1^m|^2\sigma_{s^m}^2},
\end{multline}
where $({\bf q}^m)^H=[q_1^m,\,({\bf q}_2^m)^H]= 
[(\widehat{\bf w}^m)^H{\bf a}^m,\,(\widehat{\bf w}^m)^H{\bf Q}^m]$. Assuming ``small'' estimation errors, i.e., ${\bf q}^m\approx{\bf e}_1$, similar approximation to that in \eqref{ISRinst} gives
\begin{equation}
\isr\approx \frac{1}{\sum_{m=1}^M\sigma_{s^m}^2}\sum_{m=1}^M{\tt tr}\Bigl(\C_{\z^m}{\bf q}_2^m({\bf q}_2^m)^H\Bigr).
\end{equation}
Using the equivariance property described in   Section~\ref{section:CRIB}, the CRIB is, in general, obtained through
\begin{equation}\label{CRiBM}
 {{\rm E}}\left[\isr\right] \geq \frac{1}{\sum_{m=1}^M\sigma_{s^m}^2}
 {\tt tr}\left(\sum_{m=1}^M\C_{\z^m}\crb\left(\h^m\right)|_{\substack{\h^m=0\\{\bf g}^m=0}}\right).
\end{equation}

\subsection{Block-wise ICE}
To extract the SOI from each block of data \eqref{newmodel}, the ICE approach can be used. Then, the mixing and separating vectors are estimated as parameters that are independent of  the other blocks. We will refer to this approach as to BICE.

Assuming that the background is circular Gaussian, the CRIB for BICE follows from the results of Section~\ref{section:gausscase}. By putting \eqref{CRB} into \eqref{CRiBM} and using the fact all data are independently distributed, the CRIB is given by
\begin{equation}\label{CRIBMICE}
{\rm E}[\isr]\geq \frac{1}{N_b}\frac{d-1}{\sum_{m=1}^M \sigma_{s^m}^2}\sum_{m=1}^M\frac{\sigma_{s^m}^2}{\kappa_{s_m}\sigma_{s^m}^2-1}.
\end{equation}
It is worth comparing this bound with CRIBs derived for the CMV and CSV models given by \eqref{CMV} and \eqref{CSV}, respectively, which is the subject of the following subsections. 

\subsection{Constant Mixing Vector}
In the CMV model, ${\bf a}$ is constant over $M$ blocks while the separating vector can be varying from block to block. Therefore, there are $M(d-1)+d$ free parameters. The scaling ambiguity can be resolved by putting $\gamma=1$, which is the first element of ${\bf a}$, so there are finally $(M+1)(d-1)$ free (complex-valued) parameters in the mixing model represented by parameter vectors $\g$ and $\h=[\h^1;\dots;\h^M]$.

From \eqref{likelihood}, it follows that the log-likelihood function for one sample data of the $m$th block is given by 
\begin{equation}\label{loglikelihood of mth block}
    \mathcal{L}^m(\x^m|\g,\h) = \log p_{s^m}\left(({\bf w}^m)^H\x^m\right) + \log p_{\z^m}({\bf B}\x^m),
\end{equation} 
Since the data are i.i.d. inside each block and independently distributed among the blocks, the likelihood function of the entire batch of data is equal to 
\begin{equation}
N_b\sum_{m=1}^M\mathcal{L}^m(\x^m|\g,\h). 
\end{equation}
The derivatives of \eqref{loglikelihood of mth block} are  computed similarly to \eqref{dg} and \eqref{dh}, that is,
\begin{align} \label{dgm}
\bnabla_\g^m = \frac{\p \mathcal{L}^m(\x^m|\g,\h)}{\p \g^*}\Big|_{\h={\bf 0}} &= -\boldsymbol{\psi}_{\z^m}{s^m}^* , \\ \label{dhm}
\bnabla_\h^{m,n} = \frac{\p \mathcal{L}^m(\x^m|\g,\h)}{\p {\h^n}^*}\Big|_{\h={\bf 0}} &= \delta_{n,m}\psi_{s^m}^*\z^m ,
\end{align} 
where $\boldsymbol{\psi}_{\z^m} = -\frac{\p \ln p_{\z^m}}{\p \z^*}$, $\psi_{s^m} = -\frac{\p \ln p_{s^m}}{\p s^*}$, 
and $\delta_{n,m}$ stands for the Kronecker delta. 

Now, the FIM of data from all blocks is a square matrix of dimension $2(m+1)(d-1)$ consisting of $(m+1)\times(m+1)$ blocks each of dimension $(d-1)\times(d-1)$. Let $\bnabla^m = [\bnabla_\g^m; \bnabla_\h^{m,1}; \dots; \bnabla_\h^{m,M}]$. The FIM has the structure 
\begin{equation}\label{blockFIM}
    \mathcal{J}(\g,\h) = N_b\sum_{m=1}^M \mathcal{J}^m(\g,\h) =
N_b\begin{pmatrix}
 \textbf{F} & \textbf{P} \\
 {\textbf{P}}^* & {\textbf{F}}^*
\end{pmatrix}, 
\end{equation} 
where 
\begin{equation} 
\mathcal{J}^m(\g,\h) = 
\begin{pmatrix}
 \textbf{F}^m & \textbf{P}^m \\
 {\textbf{P}^m}^* & {\textbf{F}^m}^*
\end{pmatrix} 
\end{equation} 
is the FIM for one sample of the $m$th block, and 
\begin{equation}\label{partitioned}
    \textbf{F}^m = {\rm E}\left[\bnabla^m (\bnabla^m)^H \right], \qquad 
    \textbf{P}^m = {\rm E}\left[\bnabla^m (\bnabla^m)^T \right] . 
\end{equation}
The structures of $\textbf{F}^m$ and $\textbf{P}^m$ are described in details in Appendix~B.
From them it follows that the blocks of \eqref{blockFIM} are, respectively, equal to 
\begin{equation}\label{blockFIMclosedForm}
    \textbf{F} =  \begin{pmatrix}
    \sum_{m=1}^M\boldsymbol{\kappa}_{\z^m}\sigma_{s^m}^2    &     -{\bf I}_{d-1}      &      \dots    &       -{\bf I}_{d-1} \\
    -{\bf I}_{d-1}              & {\kappa}_s^1\C_{\z^1} &    &  \bigzero          \\
    \vdots         &   \bigzero &      \ddots   &          \\
    -{\bf I}_{d-1}            &        &       & {\kappa}_s^M\C_{\z^M}     
    \end{pmatrix},
\end{equation}    
and $\textbf{P}$ is diagonal
\begin{equation}
\textbf{P} = \begin{pmatrix}
       \sum_{m=1}^M 
     {\rm E}[\boldsymbol{\psi}_{\z^m}^2]{\rm E}[({s^m}^*)^2]
     \qquad \qquad \\ 
   \qquad {\rm E}[(\boldsymbol{\psi}_{s^1}^*)^2]{\rm E}[({\z^1})^2] \qquad \\ 
  \qquad \ddots, \qquad\\
    \qquad \qquad\qquad{\rm E}[(\boldsymbol{\psi}_{s^M}^*)^2]{\rm E}[({\z^M})^2] \end{pmatrix},
\end{equation}
where 
$
\kappa_{s^m} = {\rm E}[|\boldsymbol{\psi}_{s^m}|^2]$, $\bkappa_{\z^m} = {\rm E}[\boldsymbol{\psi}_{\z^m}\boldsymbol{\psi}_{\z^m}^H]$,
 $\sigma_{s^m}^2 = {\rm E}[|s^m|^2]$, $\C_{\z^m} = {\rm E}[\z^m (\z^m)^H]$. 

For the sake of simplicity, we will consider only the special case when the  background is circular Gaussian. Then, similar simplifications to those in Section~\ref{section:gausscase} hold, ${\bf P}={\bf 0}$, $\boldsymbol{\kappa}_{\z^m} = \C_{\z^m}^{-1}$, and the  block of $\mathcal{J}^{-1}$ corresponding to $\h^m$ is
\begin{multline}\label{CRLBCMV}
    \crb (\h^m)|_{\h=0} = \frac{1}{N_b}\Bigg\{\frac{1}{\kappa_{s^m}}\C_{\z^m}^{-1} + \\ \frac{1}{\kappa_{s^m}}\C_{\z^m}^{-1} \left(\sum_{i=1}^M \frac{\sigma_{s^i}^2 \kappa_{s^i}-1}{\kappa_{s^i}}\C_{\z^i}^{-1} \right)^{-1}\frac{1}{\kappa_{s^m}}\C_{\z^m}^{-1}\Bigg\}.
\end{multline} 
By combining \eqref{CRiBM} and \eqref{CRLBCMV}, the CRIB says that
\begin{multline}\label{CRIBCMV}
    {{\rm E}}\left[\isr\right] \geq \frac{1}{N_b\sum_{m=1}^M\sigma_{s^m}^2}\sum_{m=1}^M \frac{1}{\kappa_{s^m}}\times \\ \tr \left(  \textbf{I}_{d-1}+
    \left(\sum_{i=1}^M \frac{\okappa_{s^i}-1}{\kappa_{s^i}}\C_{\z^i}^{-1} \right)^{-1}\frac{1}{\kappa_{s^m}   }\C_{\z^m}^{-1} \right).
\end{multline}

\subsection{Constant Separating Vector}
In the CSV mixing model \eqref{CSV}, ${\bf w}$ is constant over the blocks while the mixing vector can be varying. Therefore, the scaling ambiguity can be resolved by putting $\beta = 1$ while considering $\gamma^1,\dots,\gamma^M$ as dependent variables, where by \eqref{betagamma} it follows that $\gamma^m = 1-\h^H\g^m$. The free parameter vectors of the model are $\g = [\g^1;\dots;\g^M]$ and $\h$.  

Using \eqref{p_x}, the log-likelihood function for one sample of the $m$th block is
\begin{multline}\label{loglikeCSV}
     \mathcal{L}^m(\x^m|\g,\h) = \log p_{s^m}\left({\bf w}^H\x^m\right) + \log p_\z^m({\bf B}^m\x^m) + \\ + 2(d-2)\log\left|1-\h^H\g^m  \right|,
\end{multline} 
where we use the identity $\det\left(\W_{\rm ICE}\right) = (-1)^{d-1}(1-\h^H\g^m)^{d-2}$.

The structure of the FIM is the same as for the CMV model, described by \eqref{blockFIM}-\eqref{partitioned}. The blocks of \eqref{blockFIM} are given by 
\begin{equation}\label{blockFIMclosedFormCSV}
    \textbf{F} =  \begin{pmatrix}
     \boldsymbol{\kappa}_{\z^1}\sigma_{s^1}^2 &       &\bigzero    &  -{\bf I}_{d-1}     \\
      \bigzero &      \ddots   &   & \vdots       \\
    &   & \boldsymbol{\kappa}_{\z^M}\sigma_{s^M}^2 &  -{\bf I}_{d-1}  \\
-{\bf I}_{d-1}   &    \dots    & -{\bf I}_{d-1}   &            \sum_{m=1}^M{\kappa}_s^m\C_{\z^m}
    \end{pmatrix}, 
\end{equation}
and $\textbf{P}$ is diagonal
\begin{equation}\label{PCSV}
      \textbf{P} = \begin{pmatrix}
   {\rm E}[(\boldsymbol{\psi}_{\z^1})^2]{\rm E}[({s^1}^*)^2] \qquad \qquad \qquad\qquad \\ 
   \qquad \ddots, \qquad\\
    \qquad \qquad {\rm E}[(\boldsymbol{\psi}_{\z^M})^2]{\rm E}[({s^M}^*)^2] \\
\qquad \qquad     \sum_{m=1}^M  {\rm E}[(\boldsymbol{\psi}_{s^m}^*)^2]{\rm E}[({\z^m})^2]\\ 
    \end{pmatrix}.
\end{equation}

%The estimated parameter vector $\h$ corresponds the sub-matrix $\sum_{m=1}^M\boldsymbol{\kappa}_s^m\C_{\z^m}$. 
Here, we also consider only the special case that the background is circular background, for which ${\bf P}={\bf 0}$, $\boldsymbol{\kappa}_{\z^m} = \C_{\z^m}^{-1}$. Then, $\crb (\h)|_{\h=\g=0}$ is obtained as the block of the inverse matrix of FIM corresponding to the lower right-corner block of ${\bf F}$, which gives
\begin{equation}
    \crb (\h)|_{\h=\g=0} =  \frac{1}{N_b}\left(\sum_{m=1}^M{\kappa}_s^m\C_{\z^m} - \frac{1}{\sigma_{s^m}^2} \C_{\z^m} \right)^{-1}.
\end{equation} 
By putting this result in \eqref{CRiBM}, the CRIB says that
\begin{equation}\label{CRIBCSV}
    {{\rm E}}\left[\isr\right] \geq \frac{1}{N_b\sum_{m=1}^M\sigma_{s^m}^2} \tr \left( \left(\sum_{m=1}^M \frac{{\okappa}_{s^m}-1} {\sigma_{s^m}^2} \C_{\z^m} \right)^{-1}
    \sum_{m=1}^M\C_{\z^m}
    \right).
\end{equation}

\section{Discussion}\label{sec:discussion}
The expressions in brackets in \eqref{CRIBCMV} and \eqref{CRIBCSV} subject to the matrix inverse operation are non-negative combinations of positive definite matrices ($\C_{\z_m}^{-1}$ or $\C_{\z_m}$). It follows that the sums are also positive definite unless all coefficients of the linear combinations are zero. The latter case appears only if $\okappa_{s^m}=1$ for all $m$, that is, when the SOI is Gaussian on all blocks. Otherwise, the CRIBs are finite.

In the following, we discuss several special cases in order to compare the derived bounds. 

\subsection{Only one block}
When $M=1$, the piecewise determined models coincide with the standard ICE model. The reader can easily verify that, for this particular case, the bounds given by \eqref{CRIB}, \eqref{CRIBMICE}, \eqref{CRIBCMV} and \eqref{CRIBCSV} coincide as well. 

In further discussions, we will assume that $M>1$.

\subsection{An i.i.d. SOI} When the SOI has the same pdf (and also variance) in all blocks, we can denote $\kappa_s^m=\kappa_s$ and $\sigma^2_{s^m}=\sigma^2_s$ since these statistics become independent of $m$. Then, the CRIBs \eqref{CRIBMICE}, \eqref{CRIBCMV} and \eqref{CRIBCSV} can be, respectively, written in the form
\begin{align}
\text{BICE:}\quad&{\rm E}[\isr]\geq\frac{M}{N}\frac{d-1}{\okappa_s-1},\label{BICE_case2}\\
\text{CMV:}\quad&{\rm E}[\isr]\geq \frac{d-1}{N}\left(\frac{1}{\okappa_s-1}+\frac{M-1}{\okappa_s}\right),\label{CMV_case2}\\
\text{CSV:}\quad&{\rm E}[\isr]\geq \frac{1}{N}\frac{d-1}{\okappa_s-1}\label{CSV_case2},
\end{align}
A necessary condition for the identifiability of these models is that $\okappa_s>1$, which means that the SOI must have non-Gaussian pdf. The CRIB for BICE is always higher than those for CSV and CMV, which is caused by the higher complexity of BICE. CSV and CMV take advantage of the joint parameters.   

\subsection{SOI with varying variance}
Let the variance of the SOI be changing from block to block while the normalized pdf of the SOI be constant. It means that $\sigma_{s^m}^2$ depends on $m$ while $\kappa_{s^m}\sigma_{s^m}^2=\okappa_s$ is constant over the blocks. Then, the CRIBs can be written as
\begin{align}
\text{BICE:}\quad&{\rm E}[\isr]\geq\frac{M}{N}\frac{d-1}{\okappa_s-1},\label{BICE_case3}\\
\text{CMV:}\quad& {\rm E}[\isr]\geq\frac{M(d-1)}{N\okappa_s} + \frac{M}{N\okappa_s(\okappa_s-1)}
\label{CMV_case3} T_{\rm CMV},\\
\text{CSV:}\quad&{\rm E}[\isr]\geq \frac{M}{N(\okappa_s-1)}T_{\rm CSV}
\label{CSV_case3},
\end{align}
where 
\begin{align}
    T_{\rm CMV} &=\tr \left(\sum_{m=1}^M  \frac{\sigma_{s^m}^2}{\sum_{j=1}^M\sigma_{s^j}^2} \left(\sum_{i=1}^M {\bf S}_i  \right)^{-1} {\bf S}_m \right),\\
    T_{\rm CSV} &=\tr \left( \frac{1}{\sum_{j=1}^M\sigma_{s^j}^2} \left(\sum_{i=1}^M \frac{1}{\sigma_{s^i}^2}\C_{\z^m} \right)^{-1}
    \sum_{m=1}^M\C_{\z^m},
    \right)
\end{align}
and ${\bf S}_m=\sigma_{s^m}^2\C_{\z^m}^{-1}$.

The bound given by \eqref{BICE_case3} coincides with \eqref{BICE_case2}, which means that the dynamic envelop of the SOI does not have any influence on the achievable performance when ICE is independently applied to each block. By comparing \eqref{CMV_case3} with \eqref{CMV_case2} and \eqref{CSV_case3} with \eqref{CSV_case2}, we obtain more interesting results. 

To analyze, the following inequalities are needed (see Appendix~B for proofs):
\begin{align}
    \frac{d-1}{M}\leq T_{\rm CMV}&\leq d-1\label{CMV_case3_inequality} \\
    T_{\rm CSV} &\leq \frac{d-1}{M}. 
\end{align}
It follows that the bound \eqref{CSV_case3} is always lower than the one given by \eqref{CSV_case2}, moreover, by the proof it follows that the equality holds if and only if $\sigma_{s^m}^2$ is constant. It means that the non-stationarity of the SOI improves the blind extraction under the CSV model. This is not that surprising because similar conclusions follow from Cram\'er-Rao induced bounds for the standard BSS models that involve signals' non-stationarity, where more dynamical signals improve the achievable separation accuracy; see, e.g., \cite{pham2001,Koldovsky2009}.

However, by putting the lower and upper limits in \eqref{CMV_case3_inequality}, i.e. $T_{\rm CMV}=(d-1)/2$ and $T_{\rm CMV}=d-1$, into \eqref{CMV_case3}, the bound coincides with \eqref{CMV_case2} and \eqref{BICE_case3}, respectively. On one hand, it means that the achievable ISR by CMV is never worse than that by BICE. On the other hand, \eqref{CMV_case3} coincides with \eqref{CSV_case2} only if $\sigma_{s^m}^2$ is constant and is always greater than \eqref{CSV_case2} when $\sigma_{s^m}^2$ is changing. It means that the nonstationarity of the SOI is worsening the extraction accuracy under the CMV model!

\subsection{All but one blocks of SOI are circular Gaussian }
When the SOI has the circular Gaussian pdf on the $k$th block, then $\okappa_{s^k}=1$. Hence, the CRIB \eqref{CRIBMICE} does not exist when there is a block where the SOI is circular Gaussian. However, CRIBs \eqref{CRIBCMV} and \eqref{CRIBCSV} exist if the SOI is non-Gaussian or non-circular on one block at least. In the special case when all block of SOI but the $k$th block have circular Gaussian pdf, the CRIBs \eqref{CRIBCMV} and \eqref{CRIBCSV} are   
\begin{align}
    \text{CMV:} \quad {\rm E}[\isr]\geq & \frac{1}{N_b}\frac{1}{\sum_{m=1}^M\sigma_{s^m}^2}\times \\ \nonumber & \tr\left( \sum_{m=1}^M \frac{1}{\kappa_{s^m}}{\bf I}_{d-1}+\frac{\kappa_{s^k}}{\okappa_{s^k}-1}{\bf C}_{\z^k}\sum_{m=1}^M \frac{1}{\kappa_{s^m}}{\bf C}_{\z^m}^{-1} \right) ,  \\
    \text{CSV:}\quad {\rm E}[\isr]\geq & \frac{1}{N_b}\frac{1}{\sum_{m=1}^M \sigma_{s^m}^2}\frac{\sigma_{s^k}^2}{\okappa_{s^k}-1}  \tr\left({\bf C}_{\z^k}^{-1} \sum_{m=1}^M {\bf C}_{\z^m}\right). 
\end{align}    Thus, for the identifiability of CVM and CSV models is only sufficient that the SOI is not circular Gaussian on at least one block, which is a significant enhacement in comparison to BICE model.   

\subsection{Gaussian SOI and vanishing background}
Let us assume all blocks of SOI circular Gaussian and a vanishing background on the $k$th block given by 
\begin{equation}
    {\bf C}_{z^k} = \frac{\okappa_{s^k}-1}{\kappa_{s^k}}{\bf T},
\end{equation} where $T$ is a PDF matrix. In this special case the CRIBs \eqref{CRIBMICE} and \eqref{CRIBCSV} do not exist, but CRIB \eqref{CRIBCMV} is given by 
\begin{equation}
    \text{CMV:} \quad{\rm E}[\isr]\geq \frac{d-1}{N_b}
    \frac{1}{N_b}\frac{1}{\sum_{m=1}^M \sigma_{s^m}^2} \tr\left({\bf T}\sum_{m=1}^M (\sigma_{s^m}^2)^2 {\bf C}_{\z^m}^{-1}  \right).
\end{equation}

\section{Simulations}\label{secExp}
In simulations, we compare the theoretical bounds with empirical mean ISR achieved by selected ICA/ICE algorithms. Here, we have to cope with the permutation ambiguity, which means that a given algorithm need not converge to the desired SOI. In case of BSE/ICE algorithms that extract only one source, the convergence is arranged through proper initialization. For ICA methods, the SOI is identified as the separated signal with the lowest ISR.   Since the algorithms do not converged to the right SOI in some runs, the trimmed mean of ISR is computed instead of the mean, that means $10\%$ of lowest and greatest values of ISR are discarded. Replacing the mean with the trimmed mean can slightly affect the results by introducing a small bias.

\subsection{Determined Mixing Model}
\subsubsection{Gaussian Background}
Here, the CRIB given by \eqref{CRIB} assuming circular Gaussian background signals is compared with the empirical mean ISR achieved by three methods. The first, non-circular FastICA (NC-FastICA) from \cite{adaliFastICA}, is an ICA algorithm designed particularly for signals belonging to the complex Generalized Gaussian Distribution (GGD) family \cite{cGGD}, which involves also non-circular signals. The second, OGICE (Orthogonally Constrained ICE) from \cite{eusipco2017} is compared, which is an ICE algorithm derived based on maximum likelihood principle. The third, the Natural Gradient (NG) \cite{amari1996}, is a popular ICA algorithm. In OGICE, the background is modeled as circular Gaussian, therefore, this method can attain the CRIB asymptotically when provided that it is always initialized in the region of convergence to the SOI and the true score function related to its pdf is used as the internal nonlinear function.
%, i.e., within the area of convergence to the desired (target) signal. 

In trials, $d = 5$ independent complex-valued signals are generated. The target signal is drawn from the complex-valued GGD with zero mean, unit variance, shape parameter $\alpha\in(0,+\infty)$, and a circularity coefficient $\gamma \in [0,1]$. The corresponding pdf is given by \cite{loeschCRB}
\begin{equation}\label{GGDpdf}
p(s,s^*) = \frac{\alpha \rho \exp\left(-\left[\frac{\rho/2}{\gamma^2-1}\left(\gamma s^2+\gamma (s^*)^2 -2ss^*\right)\right]^\alpha\right)}{\pi \Gamma(1/\alpha)(1-\gamma^2)^{\frac{1}{2}}},
\end{equation} where $\rho = \frac{\Gamma(2/\alpha)}{\Gamma(1/\alpha)}$, and $\Gamma(\cdot)$ is the Gamma function.
The other (background) signals are circular Gaussian, which corresponds to $\alpha=1$ and $\gamma = 0$. All signals are mixed by a random mixing matrix $\A$ with elements drawn from $\mathcal{CN}(0,1)$. 
%Note that such mixture corresponds with \eqref{mixingICE} where $\C_\z%$ is determined by $\A$.

OGICE is initialized by a randomly perturbed first column of $\A$, Natural Graident is initialized by the randomly perturbed mixing matrix $\A$, while the 
initialization of NC-FastICA is random in full. In OGICE and NG, the 
nonlinearity is the same as the true score function corresponding to \eqref{GGDpdf}, that is,
\begin{equation}\label{trueSkore}
\psi(s,s^*) = \frac{2\alpha(\rho/2)^\alpha}{(\gamma^2-1)^\alpha}\left(\gamma 
s^2+\gamma (s^*)^2 -2ss^*\right)^{\alpha-1}\left(\gamma s-s^*\right).
\end{equation}
%where the parameters are set to their true values. 
It can be shown that %For the GGD pdf \eqref{GGDpdf}, $\kappa \sigma_s^2 = \okappa$ is 
\cite{loeschCRB}
 \begin{equation}
 \okappa = {\rm E}\left[\left|\psi(s) \right|^2 \right] = \frac{\alpha^2 
 \Gamma(2/\alpha)}{(1-\gamma^2)\Gamma^2(1/\alpha)}.
\end{equation}
Finally, note that NC-FastICA is endowed by the nonlinearity proposed in \cite{adaliFastICA}, the best accuracy is achieved when using kurtosis.  

\begin{figure}[!ht]
 	\centering
  \includegraphics[width=8cm]{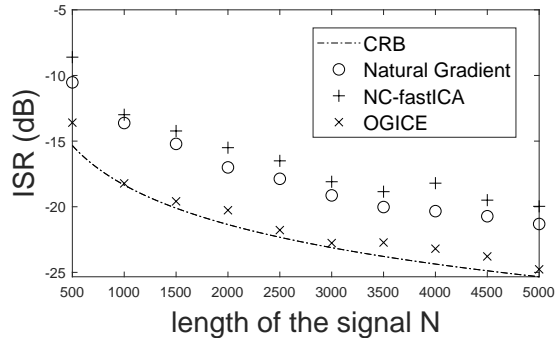}
%   LesserNvar
  \caption{Average ISR for $d=5$, $\alpha=2$, and varying 
  $N$.}\label{naN}
  \end{figure}

\begin{figure}[!ht]
 	 	\centering
  \includegraphics[width=8cm]{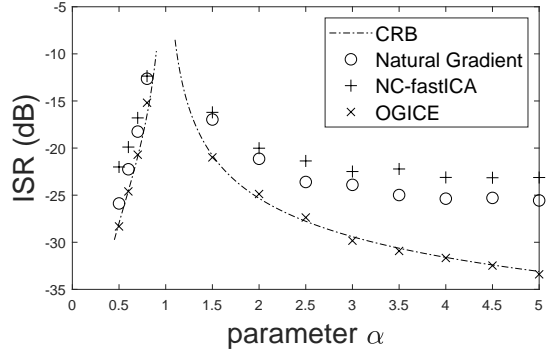}
%   LesserAlphaSub
  \caption{Average ISR for $d=5$, $N = 2500$ and varying 
  $\alpha$.}\label{naAlpha}
  \end{figure}

\begin{figure}[!ht]
 	 	\centering
  \includegraphics[width=8cm]{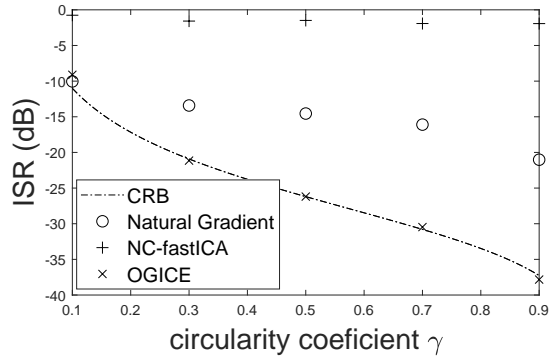}
%   lesserCircGauss
  \caption{Average ISR for $d=5$, $N = 2500$, $\alpha = 1$ and varying circularity coefficient
  $\gamma$.}\label{naCircu}
  \end{figure}

Figs.~\ref{naN}, \ref{naAlpha}
% , \ref{naAlphaSup}
and \ref{naCircu}, show 
average ISR achieved by the algorithms in $100$ trials, respectively, for 
varying $N$, %(fixed $\alpha$, $\gamma$
%, i.e. complex Laplacian random variable
%)
$\alpha$, % (fixed $N$, $\gamma$) 
and $\gamma$. % (fixed $N$, $\alpha$). 
%All variables in \ref{naN},  \ref{naAlpha} and \ref{naAlphaSup} are generated as complex circular.
The average ISRs achieved by OGICE are very close to the bound \eqref{CRIB}, 
where,   
The performance of NC-FastICA appears to be slightly limited in comparison to the NG due the versatility of 
the nonlinearity.
% , which is equal to the true score function of the SOI only 
% when
%In Fig. \ref{naAlpha}, kurtosis was used as the nonlinearity in 
%NC-FastICA, which corresponds to the true score function only when 
% $\alpha=\frac{3}{2}$. 

In Fig.~\ref{naAlpha}, 
% and \ref{naAlphaSup}
the ISR for sub-gaussian 
($\alpha>1$) and super-Gaussian ($\alpha<1$) SOI is shown.   
%for GGD with 
%$\alpha>1$, is shown. Fig. \ref{naAlphaSup} shows the accuracy for 
%super-Gaussian variables ($\alpha<1$). 
For $\alpha=1$, all signals, including the SOI, are circular Gaussian, in which 
case the mixing coefficients are not identifiable. Therefore, the ISRs approach 
$0$~dB, which means no separation. 
%With growing $N$, the bound and the ISRs are decreasing. 

In Fig.~\ref{naCircu}, the non-circular Gaussian SOI with varying circularity 
is considered. Note, that NC-FastICA is designed to be robust to circularity changes, but do not benefit from non-circularity. Thus, it does not show any dependence on 
$\gamma$, but cannot extract circular Gaussian SOI, which agrees with \cite{adaliFastICA}. The ISR achieved by OGICE 
approaches the CRIB, which confirms that a non-circular Gaussian signal can be 
extracted from the other Gaussian signals when their circularity coefficient is 
different. This condition becomes violated as $\gamma$ approaches $0$, which 
corresponds with the decaying ISR.

%   \begin{figure}[!ht]
%  	 	\centering
%   \includegraphics[width=8cm]{LesserAlphaSup}
%   \caption{Average ISR for $d=5$, $N = 2500$ and varying 
%   $\alpha$.}\label{naAlphaSup}
%   \end{figure}
  
\subsubsection{Non-Gaussian Background}
As shown in Section~\ref{sec:circularsources}, there is a coincidence between the CRIBs for ICA and ICE when, in ICE, the non-Gaussianity of background is taken into account. %In ICE, only the independence of the SOI and the background is needed. 
In this section, we simulate the case mentioned at the end of that section, that is, when background signals are dependent (a transformation decomposing them into independent components as assumed in ICA need not to exist).  The theoretical CRIB for this simulation is given by \eqref{CRLB_ngICE}.
 
In a trial, $d=4$ real-valued signals are generated. The background is 
 drawn according to the joint pdf given by 
\begin{equation}\label{joined} 
p(z^1,\dots,z^{d-1}) \propto \exp\left(-\left({\lambda \sum_{i=1}^{d-1} |{z^i}|^2}\right)^{\alpha}\right)
\end{equation}        where %$D_1$ is a normalization constant, 
$\lambda>0$, and $\alpha \neq 1$ (for $\alpha = 1$, the pdf is Gaussian). To scale the marginal pdfs of background signals to the unit variance, we put $
\lambda = \frac{\Gamma\left(\frac{5}{2\alpha} \right)}{3\Gamma\left(\frac{3}{2\alpha} \right)}
$. Then, it holds that 
\begin{equation}
    (\bkappa_{\z})_{kk}= %\kappa_{\bf z}^k = {\rm E}\left[{\psi_{\rm IVE}^k}^2 \right] =
    \frac{4}{3}\lambda \alpha^2\frac{\Gamma(2+\frac{1}{2\alpha})}{\Gamma(\frac{3}{2\alpha})}. 
\end{equation}
The SOI is drawn from the GGD with zero mean, unit variance and a shape parameter $\tilde{\alpha}$, where $\tilde{\alpha} = \alpha + 1$. Thus, for $\alpha=1$, all signals in the mixture are Gaussian. 

We compare three algorithms with the CRIB given by \eqref{CRLB_ngICE}: OGICE \cite{eusipco2017}, EFICA \cite{koldovsky2006}, and NG-OGICE \cite{koldovsky2018a}. OGICE is designed for ICE with Gaussian background, where the CRIB is given by \eqref{CRIB} (which we show as well for the sake of completeness). EFICA is an asymptotically efficient ICA algorithm provided that all original signals are drawn from GGD. %It is a method for ICA whose bound is given by \eqref{CRLB}.
NG-OGICE is an ICE method considering the non-Gaussianity of background, in which the true multivariate score function of background must be known. %The bound for non-Gaussian ICE \eqref{CRLB_ngICE} considers possibly dependent signals in the background. 

\begin{figure}[!ht]
 	 	\centering
  \includegraphics[width=8cm]{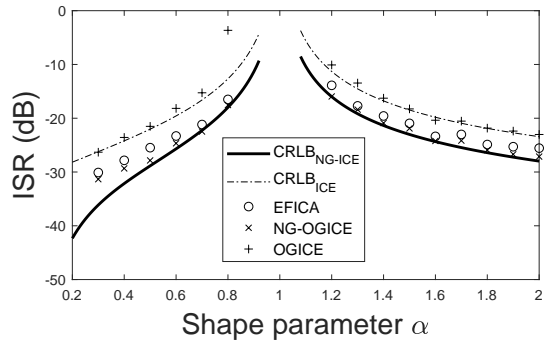}
  \caption{Average ISR for non-Gaussian background when pdfs of all signals are varying with respect to $\alpha$.}\label{fig:ng_background}
  \end{figure}

In Fig.~\ref{fig:ng_background}, the ISRs averaged over 100 trials achieved by OGICE, EFICA and NG-OGICE are compared. The bound \eqref{CRLB_ngICE} is denoted by CRIB$_{\rm NG-ICE}$ and the one for the Gaussian background \eqref{CRIB} is denoted by CRIB$_{\rm ICE}$. The results show that the mean ISRs by OGICE are close to the bound given by \eqref{CRIB} (which is in a good agreement with the results of asymptotic performance analyses \eqref{CRLBICE} \cite{pham2006}). The results by EFICA and NG-OGICE are closer to \eqref{CRLB_ngICE}. NG-OGICE is even slightly more accurate than EFICA, which is caused by a more accurate modeling of the background's pdf.  

For $\alpha=1$,  all signals are Gaussian, which means that the SOI cannot be separated from the background. With increasing non-Gaussianity of the mixture, which means increasing distance from $\alpha=1$, the separation accuracy gets better. 

\subsection{Piecewise Determined Mixtures with Circular Gaussian Background}
To validate the bounds for CMV and CSV, both are compared with empirical results achieved by block-wise versions of OGICE introduced in \cite{koldovsky2019icassp}. The methods will be jointly referred to as BOGICE (in \cite{koldovsky2019icassp}, BOGICE$_{\bf a}$ is the variant for CMV while BOGICE$_{\bf w}$ is for CSV). It should be noted that no other methods for CMV/CSV currently exist in the literature to our best knowledge. 

In experiments here, we consider two statistical models of signals: The SOI is i.i.d. non-Gaussian over all blocks and i.i.d. SOI within blocks with the same distribution but varying variance over blocks. The background is assumed circular Gaussian i.i.d. with unit variance in all blocks in both cases.

In trials, $d = 5$ independent complex-valued signals are generated. The SOI is drawn from a circular complex GGD with zero mean, unit variance, $\alpha=2$. The other signals are circular Gaussian, which corresponds to $\alpha=1$. The nonlinearity is given by the true score function. $M$ blocks of the same length are considered. Each block is mixed by a random mixing matrix. The mixing matrices obey the mixing models CMV or CSV, respectively. 

The empirical ISRs achieved by BOGICE are compared with the CRIB corresponding to the mixing model used in the given simulation. To compare, we always show also the hypothetical CRIB achieved when the alternative mixing model was considered with the same statistical properties of the SOI.

\subsubsection{An i.i.d. SOI}
\begin{figure}[!ht]
 	 	\centering
  \includegraphics[width=8cm]{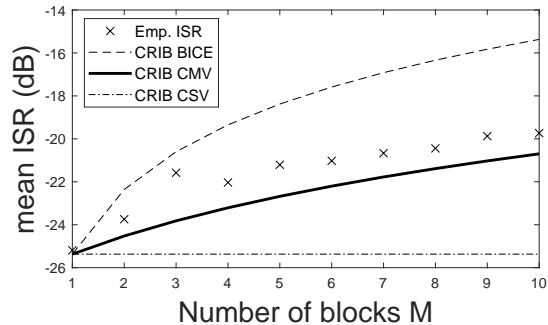}
  \caption{Average ISR for CMV mixing model when $d=5$, $N = 5040$, and varying number of blocks $M$.}\label{fig:cmvMvar}
  \end{figure}

\begin{figure}[!ht]
 	 	\centering
  \includegraphics[width=8cm]{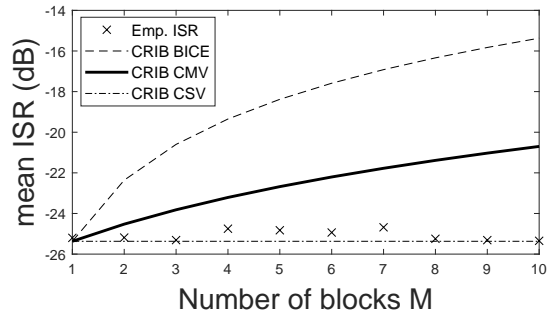}
  \caption{Average ISR for CSV mixing model when $d=5$, $N = 5040$, and varying number of blocks $M$.}\label{fig:csvMvar}
\end{figure}

Fig.~\ref{fig:cmvMvar} corresponds to the simulation considering the CMV model for varying number of blocks, that is, $M=1,2,5,10$. It shows bar chart of ISR achieved by BOGICE averaged over 500 trials and the CRIB given by \eqref{CMV_case2} (CMV) and, for comparison, also the CRIBs \eqref{BICE_case2} (BICE) and \eqref{CSV_case2} (CSV). 
Similar simulation was done with the CSV model; the results are shown in 
Fig.~\ref{fig:csvMvar} in the same fashion as in Fig.~\ref{fig:cmvMvar}.

Figs. \ref{fig:cmvMvar} and \ref{fig:csvMvar} show the coincidence between the empirical results by the variants of BOGICE and the CRIBs corresponding to the mixing model of the given simulation. The performances of the methods follow the same dependence on the number of blocks $M$ as these CRIBs. The results also show that BOGICE takes the advantage of the special mixing model CMV/CSV compared to BICE, as its mean ISR is lower that the CRIB \eqref{BICE_case2}, unless $M=1$ where all mixing models coincide. 

The CRIBs by CSV are lower than those for CMV, which agrees with the results of Section~\ref{sec:discussion}. However, with this conclusion, the fact that both models are based on different assumptions must be also taken into account.

\subsubsection{SOI with varying variance}
In this special case, the SOI with the same pdf but varying variance over blocks is assumed. In a trial, $M=5$ blocks and four different settings of SOI's variances are considered: Specifically, type~A is $\sigma_{s^m}^2 = 1$ for $m=1, \dots,5$, type~B corresponds to $\sigma_{s^1}^2=\sigma_{s^2}^2=1$, $\sigma_{s^2}^2=1$, $\sigma_{s^4}^2=\sigma_{s^5}^2=3$, type~C shows $\sigma_{s^m}^2=m$, $m=1,\dots,5$ and type~D is for $\sigma_{s^m}^2=m^2$, $m=1,\dots,5$. 
\begin{figure}[!ht]
 	 	\centering
  \includegraphics[width=8cm]{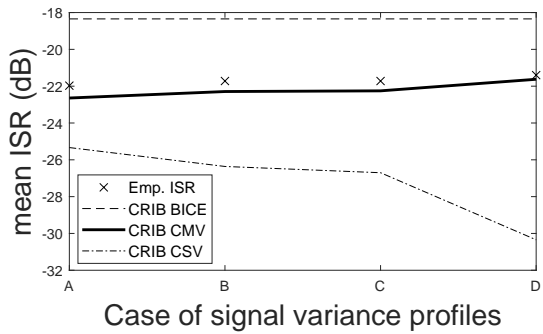}
  \caption{Average ISR for CMV mixing model when $d=5$, $N = 5000$, and varying $\sigma_{s^m}$ over blocks.}\label{fig:cmvSvar}
  \end{figure}
  
  \begin{figure}[!ht]
 	 	\centering
  \includegraphics[width=8cm]{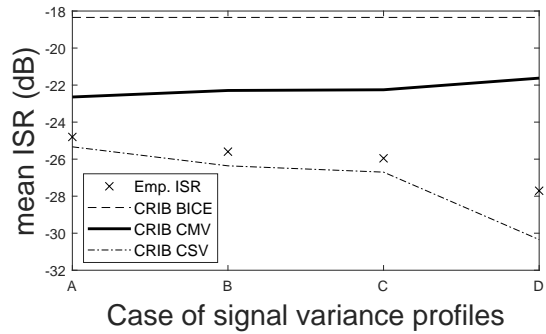}
  \caption{Average ISR for CMV mixing model when $d=5$, $N = 5000$, and varying $\sigma_{s^m}$ over blocks.}\label{fig:csvSvar}
  \end{figure} 
  
%   Type equal to $1$ in Figs. \ref{fig:cmvSvar} and \ref{fig:csvSvar} shows results when all blocks have the same variance 
  
  In \eqref{CMV_case3} and \eqref{CSV_case3} is shown that the nonstationarity of the SOI improves the separation accuracy under the CSV mixing model, but worsen the accuracy under the CMV model.
  
  The empirical results and so do the bounds show that for the increasing number of blocks $M$ the accuracy under the CMV mixing model drops, but under the CSV model it levels up. When the number of blocks is fixed and the variance of SOI differs on blocks then the results show that the higher the diversity of variance the better accuracy under CSV model, but lower under CMV model.
\section{Conclusions}\label{secConslusion}
The derived CRLB on ISR achieved by complex ICE have shown that the structured (de-)mixing matrix model with a reduced number of parameters is not restrictive. The accuracy achievable by ICE with circular Gaussian background is asymptotically the same to that one of ICA when all but one signals are circular Gaussian. 
The CRLB shows that the general lower bound for ICA can be attained by the non-Gaussian ICE.

The piecewised determined model allows us to deal with dynamic mixtures thanks to its block structure. The CRLB of this model shows how the performance is limited by the number of blocks and that it can benefit from a varying variance of signals over blocks. 
Numerical simulations have confirmed the validity of the CRIBs.

\bibliographystyle{IEEEtran}

\bibliography{IEEEfull,conferences,ISI}

\end{document}